\documentclass[a4paper,12pt,reqno]{amsart}
\usepackage{amssymb}
\usepackage{amsmath}

\def\i{\,\lrcorner\,}
\def\a{\alpha}
\def\b{\beta}
\def\c{\gamma}

\def\vs{\vskip .6cm}

\def\.{\cdot}
\def\O{\Omega}
\def\n{\nabla}

\def\t{\tilde}
\def\beq{\begin{equation}}
\def\eeq{\end{equation}}
\def\bea{\begin{eqnarray*}}
\def\eea{\end{eqnarray*}}
\def\ba{\begin{array}}
\def\ea{\end{array}}

\def\f{\varphi}

\def\e{\varepsilon}

\def\r{\end{proof}}

\def\pa{\partial}
\def\dt{\pa_t}

\def\ci{{\mathcal C}^\infty}
\def\CP{\CM{\rm P}}

\def \RM{\mathbb{R}}

\def \CM{\mathbb{C}}

\def\d{{\delta}}

\def\Ric{\mathrm{Ric}}

\def\U{\rm{U}}
\def\SU{\rm{SU}}
\def\OO{\rm{O}}

\newtheorem{ede}{Definition}[section]

\newtheorem{epr}[ede]{Proposition}

\newtheorem{ath}[ede]{Theorem}

\newtheorem{elem}[ede]{Lemma}

\newtheorem{ecor}[ede]{Corollary}

\title{Conformally Einstein Products and Nearly K\"ahler Manifolds}
\author{Andrei Moroianu}
\author{Liviu Ornea}
\thanks{ We acknowledge several valuable suggestions from the referee
  which helped us to improve the content of the paper. The second
  author thanks the Centre de Math\'ematiques de l'Ecole
  Polytechnique for hospitality during the preparation of this
  work. He was also partially supported by grant
  2-CEx-06-11-22-/25.07.2006.} 

\address{Centre de Math{\'e}mathiques, Ecole Polytechnique, 91128
  Palaiseau Cedex, France} 
\email{am@math.polytechnique.fr}
\address{Univ. of Bucharest, Faculty of Mathematics,
14 Academiei str, 70109 Bucharest, Romania}
\email{lornea@gta.math.unibuc.ro, Liviu.Ornea@imar.ro}

\begin{document}

\begin{abstract}
In the first part of this note we study compact Riemannian
manifolds $(M,g)$ whose Riemannian product with $\RM$ is conformally
Einstein. We then consider 6-dimensional almost Hermitian
manifolds of type $W_1+W_4$ in the Gray-Hervella classification
admitting a parallel vector field and show that (under some mild
assumption) they are obtained as Riemannian cylinders over compact
Sasaki-Einstein $5$-dimensional manifolds.

\vs


\end{abstract}

\maketitle

\section{Introduction}

The study of conformally Einstein metrics goes back to Brinkmann
who determined in \cite{br1} the necessary and sufficient conditions
for a Riemannian manifold to be mapped conformally on an Einstein
manifold, and considered in \cite{br2} the special case of conformal mappings
between Einstein manifolds.  

More recently, Listing \cite{li1, li2} and Gover and Nurowski
\cite{gn} have found tensorial obstructions for (semi-)Riemannian
metrics to be 
conformally Einstein under some non-degeneracy hypothesis for the Weyl
tensor.

Motivated by a problem coming from almost Hermitian geometry, we study 
conformally Einstein metrics from a different point of
view. More precisely, we look for conformally Einstein metrics of
product type $g+dt^2$ on cylinders $M\times \RM$. In Theorem \ref{prod} we
classify all such metrics in the positive scalar curvature case,
assuming that $M$ is compact. We show that $(M,g)$ has to be Einstein
with positive scalar curvature and, moreover, that the
conformal Einstein factor on $M\times \RM$ can be explicitly
determined and only depends on the $\RM$-coordinate.

Although the problem is local, we had to assume compactness in order
to solve completely the system of equations it leads to. We do not
know whether the compactness hypothesis can be removed.   

In the second part of this paper we turn our attention to nearly
K\"ahler geometry, a subject which appears to be very important in
contemporary theoretical physics (cf., for example,  \cite{fi},
\cite{ve}). After  recalling some basic facts about nearly K\"ahler
manifolds in Section \ref{3} and reviewing the Gray-Hervella
classification in Section \ref{4}, we address the question of the
existence of (6-dimensional) nearly K\"ahler conformal structures 
on cylinders over compact manifolds. The link between this and
the conformal Einstein problem is provided by the fact that nearly
K\"ahler manifolds are automatically Einstein in dimension 6.

Our main result roughly says that a Riemannian cylinder
$N^5\times\RM$ is conformally nearly K\"ahler if and only if the basis
is Sasaki-Einstein (see Theorem \ref{cylnk} for a precise statement).

\section{Conformally Einstein products}

This section is devoted to the study of conformally Einstein metrics
of cylindrical type. Notice that the special case where the conformal
factor only depends on the coordinate of the generator, corresponds to
warped products with one-dimensional basis and was studied in
\cite[9.109]{bes} and \cite[Lemma 13]{ku}.

\begin{ath}\label{prod} Let $(M^n,g)$ be a compact Riemannian
  manifold not isometric to a round sphere. Suppose that the
  Riemannian cylinder over $M$ is 
  conformally Einstein with positive scalar curvature, that is, there
  exists a smooth function $f$ on $M\times\RM$, such that the Ricci
  tensor of the Riemannian manifold
  $(M\times\RM,e ^{2f}(g+dt^2))$ is a positive multiple of the
  metric. Then $(M,g)$
  is an Einstein manifold with positive scalar curvature and the
  conformal factor satisfies  
  $e ^{2f(x,t)}=\a^2\cosh^{-2}(\b t+\c)$ for some real constants $\a,\b,\c$.
\end{ath}

\begin{proof} Let us denote by $N$ and $\t N$ the Riemannian manifolds
  $(M\times\RM,(g+dt^2))$ and $(M\times\RM,e ^{2f}(g+dt^2))$ respectively.
  We view the function $f$ on $N$ as a smooth 1-parameter family of
  functions on $M$ by $f_t(x):=f(x,t)$. In this way, the exterior
  derivative $df$ satisfies $df=df_t+f'_tdt$, where $df_t$ denotes the
  derivative of $f_t$ on $M$ and $f'_t=\frac{\pa f}{\pa t}$.
  Similarly the Laplace operators of $M$ and $N$ are related by
  $\Delta ^Nf=-f''_t+\Delta ^Mf_t$. In the sequel we shall denote by
  $\pa_t:=\frac{\pa}{\pa t}$ the ``vertical" vector field on $N$ and by
  $X,\ Y,$ vector fields on $M$, identified with their canonical
  extension to $N$ commuting with $\pa_t$. The formula for the
  conformal change of the Ricci tensor (see {\em e.g.} \cite[1.59]{bes}),
\beq\label{e}\Ric^{\t N}=\Ric^N-(n-1)
(\n^N df-df\otimes df)+(\Delta ^Nf-(n-1)|df|^2)(g+dt^2),
\eeq 
yields in particular
\beq\label{e2}\Ric^{\t
  N}(X,\dt)=-(n-1)(X(f'_t)-X(f_t)f'_t),\qquad\forall\ X\in TM.
\eeq 
Since $\t N$ is Einstein, (\ref{e2}) shows that $X(f'_t)=X(f_t)f'_t$,
which can be rewritten as $X(\dt(e ^{-f}))=0$ for every $X\in
TM$. Consequently, there exist smooth functions $a\in\ci(\RM)$ and
$b\in\ci(M)$ such that $e ^{-f(x,t)}=a(t)+b(x)$, in other words
$$f(x,t)=-\ln(a(t)+b(x)),\qquad\forall\ x\in M,\ \forall\ t\in \RM.$$
We now readily compute
\beq\label{f1}df=df_t+f'_tdt=-\frac{db+a'dt}{a+b},\qquad
f_t''=-\frac{a''(a+b)-(a')^2}{(a+b)^2},
\eeq
\beq\label{f2}\Delta ^Nf=-f''_t+\Delta
^Mf_t=\frac{a''(a+b)-(a')^2}{(a+b)^2}-\frac{\Delta 
  ^Mb}{a+b}-\frac{|db|^2}{(a+b)^2},
\eeq
and
\bea\n^Ndf&=&=-\frac{\n^N(db+a'dt)}{a+b}+\frac{1}{(a+b)^2}d(a+b)
\otimes(db+a'dt)\\ 
&=&-\frac{\n^M(db)+a''dt^2}{a+b}+\frac{1}{(a+b)^2}(db+a'dt)\otimes(db+a'dt)\\
&\stackrel{(\ref{f1})}{=}&-\frac{H(b)+a''dt^2}{a+b}+df\otimes df,
\eea
where $H(b)$ denotes the Hessian of $b$ on $M$.
Let $r$ denote the Einstein constant of $\t N$. Plugging the relations
above back into (\ref{e}) yields
\bea \frac{r(g+dt^2)}{(a+b)^2}=\Ric^{\t
  N}&=&\Ric^M+\frac{n-1}{a+b}(H(b)+a''dt^2)
\\&+&\frac{(a''-\Delta ^M
  b)(a+b)-n(a')^2-n|db|^2}{(a+b)^2}(g+dt^2),
\eea
which is equivalent to the system
$$\begin{cases}
r=(na''-\Delta^Mb)(a+b)-n(a')^2-n|db|^2\\
rg=((a''-\Delta ^M b)(a+b)-n(a')^2-n|db|^2)g+(a+b)^2\Ric^M+(n-1)(a+b)H(b).
\end{cases}$$
Subtracting  the first equation from the second one, the system becomes
\beq\label{sy}\begin{cases}
r=(na''-\Delta^Mb)(a+b)-n(a')^2-n|db|^2\\(n-1)a''g=(a+b)\Ric^M+(n-1)H(b).
\end{cases} 
\eeq
We distinguish three cases:

{\sc Case 1:} $b$ is constant on $M$. By a suitable change of
coordinates, the metric becomes a warped product and  the conclusion
could be directly derived from \cite[9.109]{bes}.  We will
nevertheless provide the direct argument. Replacing $a$ by $a-b$, we may
assume that $b=0$, so by the first equation in (\ref{sy}), $a$
satisfies the ODE 
$$a''a-(a')^2=\frac{r}{n}.$$
The general solution of this equation is 
$$a(t)=\sqrt{\frac{r}{n\b^2}}\,\cosh(\b t+\c).$$
Thus $e ^{2f(t)}=\frac{1}{a ^2}=\a^2\cosh^{-2}(\b t+\c)$, with
$\a:=\frac{n\b^2}{r}$. The second 
relation in (\ref{sy}) shows that $M$ is Einstein, with positive
Einstein constant $\b^2$. 

{\sc Case 2:} $a$ is constant on $\RM$. Again, replacing $b$ by $b-a$, we may
assume that $a=0$, so the first equation of (\ref{sy}) becomes 
$$n|db|^2+b\Delta ^Mb+r=0.$$
Integrating over $M$ yields
$$0=\int_M(n|db|^2+b\Delta ^Mb+r)dv=\int_M((n+1)|db|^2+r )dv>0,$$
showing that this case is impossible.

{\sc Case 3:} Neither $a$ nor $b$ are constant functions. We
differentiate the second relation of (\ref{sy}) twice, first with respect 
to $t$, then with respect to some arbitrary vector $X\in TM$ and obtain
\beq\label{ab}a'''X(b)=\frac1n a'X(\Delta ^M b).\eeq
Taking some $x\in M$ and $X\in T_xM$ such that $X_x(b)\ne 0$, this
relation shows that $a'''=\d a'$ for some $\d\in\RM$.
Similarly,
taking some $t\in\RM$ such that $a'(t)\ne 0$ gives some $\d'\in\RM$
such that  $X(\Delta ^M b)=\d'X(b)$ for all $X\in TM$. Plugging these
two relations back into (\ref{ab}) yields $\d'=n\d$. Summarizing, we
have
\beq\label{sys}\begin{cases}a''=\d a+\e\\ \Delta ^M b=n\d b+\e'\end{cases}\eeq
for some real constants $\e,\e'$. If $\d=0$, the second relation yields
(by integration over $M$)
$\e'=0$, so $b$ is constant, a contradiction. Thus $\d\ne
0$. Replacing $a$ by $a+\frac\e\d$ (and correspondingly replacing $b$
by $b-\frac\e\d$), we may assume $\e=0$. 
The second relation in
(\ref{sys}) also shows that $n\d$ is an eigenvalue of the Laplace
operator (corresponding to the eigenfunction $b+\frac{\e '}{n\d}$),
whence $\d>0$.
The second equation of the system
(\ref{sy}) now becomes
\beq\label{sy1}\begin{cases}
\Ric^M=(n-1)\d g\\
H(b)=-b\d g
\end{cases}\eeq
Since $b$ is non-zero, the Obata theorem (see \cite[Theorem 3]{ob})
implies that $M$ is isometric to a round sphere,
a contradiction, which shows that this case is impossible as well.
\r

We will give a concrete application of this theorem in Section
\ref{w1}, after reviewing some special classes of almost Hermitian
manifolds in the next two sections.

\section{Basics on nearly K\"ahler geometry}
\label{3}

A. Gray was led to define nearly K\"ahler manifolds (also known as
almost Tachibana spaces) by his research on weak holonomy of
$\mathrm{U}_n$-structures. An almost Hermitian manifold  $(N,h,J)$, with fundamental two-form $\Omega$ and Levi-Civita
connection $\nabla$ is called {\em nearly K\"ahler} if $\nabla
\O$ is totally skew-symmetric. 

From the viewpoint of the representations of $\mathrm{U}_n$ on the
space of tensors with the same symmetries as $\nabla \Omega$, nearly
K\"ahler manifolds appear in the class $W_1$ of the Gray-Hervella
classification (see \cite{gh} and Section \ref{4} below). It is also known
that nearly K\"ahler manifolds with integrable almost complex
structure are necessarily K\"ahler. The specific, non-trivial, case
is then the so-called \emph{strict} nearly K\"ahler, when $\nabla
J\neq 0$ at every point of $M$. 

The local structure of nearly K\"ahler manifolds was first discussed
by Gray in \cite{gr2} and was recently completely understood by
P.A. Nagy in \cite{na}: any nearly K\"ahler manifold is locally a
product of $6$-dimensional strict nearly K\"ahler manifolds, 
locally homogeneous manifolds and twistor spaces of positive
quaternionic K\"ahler manifolds. According to this result, what remains
to be studied are  strict nearly K\"ahler structures in dimension
$6$. This is, in fact, the first interesting case, since
$4$-dimensional nearly K\"ahler manifolds are automatically
K\"ahler. On the other 
hand, the dimension 6 is particularly important also because of the
following result: 
\begin{epr}\label{gra} \cite{ma} A strict $6$-dimensional nearly K\"ahler
  manifold is Einstein with positive scalar curvature. 
\end{epr}
In fact, in dimension 6, a strict nearly K\"ahler
structure is equivalent with the local
existence of a non-trivial real Killing spinor, cf. \cite{gru}.  

Nearly K\"ahler structures are closely related to
$G_2$ structures \emph{via} the cone construction (see \cite{ba} for
example): 
\begin{epr}\label{cone1} A  Riemannian manifold
  $(N^6,h)$ carries a nearly K\"ahler structure if and only if its cone
  $(N^6\times \RM^+,t^2h+dt^2)$ has holonomy contained in $G_2$. 
\end{epr}
For later use, let us also recall the following related result:
\begin{epr}\label{cone2} {\em \cite{ba}} A  Riemannian manifold 
  $(M^{2m+1},g)$ carries a Sasaki-Einstein structure if and only if its cone
  $(M\times \RM^+,t^2g+dt^2)$ has holonomy contained in $\SU_{m+1}$. 
\end{epr}
Returning to nearly K\"ahler geometry, the only known compact examples in
dimension $6$ are homogeneous: the sphere $S^6$, the 3-symmetric
space $S^3\times S^3$ and the twistor spaces, $\CP^3$ and $F(1,2)$, of
the $4$-dimensional 
self-dual Einstein manifolds $S^4$ and $\CP^2$. On the other hand,
Butruille proved recently that every compact homogeneous strictly
nearly K\"ahler manifold must be one of these (cf. \cite{jbb1}).

\section{A review of the Gray-Hervella classification}
\label{4}

Nearly K\"ahler manifolds can be understood better in terms of the
classification of almost Hermitian structures \cite{gh}.

For each almost Hermitian manifold $(N^{2m},h,J)$, with fundamental
form $\O:=h(J.,.)$, the Nijenhuis tensor, viewed as a tensor of type
$(3,0)$ {\em via} the metric, splits in two components $N=N_1+N_2$,
where $N_1$ is totally skew-symmetric and $N_2$ satisfies the Bianchi
identity.  Similarly, the covariant
derivative of $J$ with respect to the Levi-Civita connection of $h$
splits in four components under the action of the structure group
$\U_m$ (see \cite{gh}):
$$\n J= (\n J)_1+ (\n J)_2+ (\n J)_3+ (\n J)_4.$$
The first component corresponds to $N_1$, and also to the
$(3,0)+(0,3)$-part of $d\O$. The second component can be identified
with $N_2$, while the two other components correspond respectively to
the primitive part of $d\O ^{(2,1)+(1,2)}$ and to the contraction
$\O\i d\O$ which is a 1-form called the {\em Lee form}. The manifold
$N$ is called of type  $W_1+W_4$ if $(\n J)_2$ and $(\n J)_3$ vanish
identically. Similar definitions apply for every subset of subscripts
in $\{1,2,3,4\}$. For example a manifold of type $W_3+W_4$ is
Hermitian, a manifold of type $W_2$ is symplectic, and a manifold of
type $W_1$ is nearly K\"ahler.

From the definition it is more or less obvious that if the metric $h$
is replaced by a conformally equivalent metric $\t 
h:= e ^{2f}h$, the first three components of $\n J$ are invariant and
the fourth component satisfies $(\t\n J)_4=(\n J)_4+df$. Therefore, 
the Lee form of a manifold of type $W_1+W_4$ is closed (resp. exact),
if and only if the manifold is locally (resp. globally) conformal
nearly K\"ahler. In dimension 6, Butruille proved in \cite{jbb} that
all manifolds in class $W_1+W_4$ have \emph{closed} Lee form, hence
they have to be locally conformal nearly K\"ahler. But recently,
Cleyton and Ivanov proved, \cite[Lemma 8]{ci}, that every locally conformal
nearly K\"ahler structure is actually globally
conformal. Combining these results, we may state: 

\begin{ath}\label{b}{\em \cite{jbb, ci}} Let $M$ be a $6$-dimensional
almost Hermitian manifold of type $W_1+W_4$. Then its Lee form is
closed, and $M$ is either globally conformal nearly K\"ahler or locally
conformal K\"ahler.  
\end{ath}

This generalizes the well-known fact that for $m\ge 3$, every
almost Hermitian manifold of type $W_4$ is locally conformal
K\"ahler (lck).

\section{Conformally nearly K\"ahler cylinders}\label{w1}

The aim of this section is to classify all compact $5$-dimensional
Riemannian manifolds $(M,g)$ with the property that the Riemannian
cylinder $M\times\RM$ carries an almost Hermitian structure of type
$W_1+W_4$.

\begin{ath}\label{cylnk} If the Riemannian cylinder $N:=(M\times
  \RM,g+dt^2)$ over a compact 
 $5$-dimensional Riemannian manifold $(M,g)$
 carries an almost Hermitian structure of 
 type $W_1+W_4$ which is not of type $W_4$, then $(M,g)$
 is Sasaki-Einstein. Conversely, if $(M,g)$
 is Sasaki-Einstein, then its cylinder $N$ carries a structure of 
type $W_1+W_4$, besides its canonical Vaisman
structure $($locally conformal K\"ahler, \emph{i.e.} $W_4$, with
parallel Lee form, cf. \cite{ov}$)$. 
\end{ath}

\begin{proof}
Assume first that $N$ carries a structure of 
type $W_1+W_4$. By Theorem \ref{b},  $N$ 
is either globally conformal nearly K\"ahler or lck. Since we
assumed that $N$ is not lck, there exists a function $f$
on $N$ such that $(N,e ^{2f}(g+dt^2))$ is a strict nearly K\"ahler
manifold. By Proposition \ref{gra}, every such manifold in dimension $6$ is
Einstein with positive scalar curvature. We apply Theorem
\ref{prod} and obtain that either $M$ is the round sphere (which, in
particular is Sasaki-Einstein), or $e^{2f}=\a^2\cosh^{-2}(\b t+\c)$
for some real constants $\a,\b,\c$. 

Let us now consider the diffeomorphism 
$$\f:M\times\RM\to M\times (0,\pi),\qquad (x,t)\mapsto (x,2\tan^{-1}(e
^{\b t+\c})).$$ 
A straightforward computation shows that 
\beq\label{si}e ^{2f}(g+dt^2)=\frac{\a ^2}{\b^2}\f^*(\b^2\sin^2s\,g+ds^2).
\eeq
We have obtained that the so-called {\em sine-cone} (see
\cite{fimu}) of $(M,\b^2 g)$ has a nearly K\"ahler structure. The first
part of the theorem then follows from the next lemma, which can be found
(in a slightly different version) in \cite{fimu}. 

\begin{elem}\label{lem}
Up to constant re-scalings, the sine-cone
$(M\times(0,\pi),\sin^2s\,g+ds^2)$ of a 
$5$-dimensional Riemannian manifold $(M,g)$  has a nearly K\"ahler
structure if and only if $M$ is Sasaki-Einstein. 
\end{elem}
\noindent{\em Proof of the lemma.} In \cite{fimu} the authors prove the result
by an explicit calculation, using so-called {\em hypo structures} on
Sasaki-Einstein manifolds. We provide here a different argument. 
The key idea is the fact that the Riemannian product of two cone
metrics is again a cone metric, as shown by the formula
$$(t^2g+dt^2)+(s^2h+ds^2)=r^2(\sin^2\theta\,g+\cos^2\theta\,h+d\theta ^2)+dr^2,
\qquad (s,t)=(r\cos\theta,r\sin\theta).$$
In particular, taking $h=0$ (the metric of a point), shows that the
cylinder over the Riemannian cone of a metric $g$ is isometric to the
Riemannian cone over the sine-cone of $g$. By Proposition \ref{cone2}, if $M$
is Sasaki-Einstein, its Riemannian 
cone $\bar M$ has holonomy in $\SU_3$, so the cylinder $\bar M\times\RM$
has holonomy in ${\SU}_3\times\{1\}\subset G_2$. The previous remark,
together with Proposition \ref{cone1}, shows that the sine-cone of $M$ is
nearly K\"ahler. Conversely, if this holds, then the Riemannian cone
of the sine-cone has holonomy in $G_2$. Thus the holonomy of the
cylinder $\bar M\times\RM$ is a subgroup of $G_2$. But since $G_2\cap
(\OO_6\times\{1\})=\SU_3\times\{1\}\subset \OO_7$, this means that
$\bar M$ has holonomy in $\SU_3$, so 
$M$ is Sasaki-Einstein. The lemma, and the first part of the theorem
are thus proved.

Conversely, let $(M^5,g)$ be a 
Sasaki-Einstein manifold. We first notice that by Proposition \ref{cone2},
the Riemannian cone $\bar 
M$ is K\"ahler, so the cylinder $M\times \RM$, which is conformal to
$\bar M$, is lck (and even Vaisman, see \cite{ov}). 

On the other hand, Lemma \ref{lem} shows that the sine-cone of $M$
is nearly K\"ahler, and therefore the cylinder
over $M$, which by (\ref{si}) is conformal to the sine-cone, has a
structure of type $W_1+W_4$. 
\r

\noindent{\bf Remark.} As the referee pointed out, Theorem \ref{cylnk} above
has the following interesting consequence:
\begin{ecor}
If $\f$ is an isometry of a compact $5$-dimensional Riemannian manifold
$(M,g)$, the mapping torus of $\f$ does not carry any compatible
almost Hermitian structure of type $W_1+W_4$. 
\end{ecor}
\begin{proof} 
The mapping torus of $\f$ is the quotient of the Riemannian product
$(M\times\RM, g+dt^2)$ by the discrete group of isometries generated
by $(x,t)\mapsto (\f(x),t+1)$. If the quotient carries a structure of
type $W_1+W_4$, the same holds for the cylinder $M\times \RM$ (by
pull-back). From Theorem \ref{cylnk} we thus obtain that $M$ is
Sasaki-Einstein. Moreover, the almost Hermitian structure on the cylinder
can be easily made explicit, cf. \cite[Theorem 3.6]{fimu}. As it was
pointed out by S. Ivanov, this almost complex structure is not 
invariant by any translation in the $\RM$-direction, a contradiction.
\end{proof}


 \labelsep .5cm

\end{document}